\documentclass{cimart}

\usepackage{mathabx}

\DeclareMathOperator{\Aut}{Aut}
\DeclareMathOperator{\Char}{char}
\DeclareMathOperator{\SL}{SL}
\DeclareMathOperator{\tr}{tr}

\newcommand{\al}{\alpha}
\newcommand{\be}{\beta}
\newcommand{\cc}{\mathbf{c}}
\newcommand{\corr}[1]{#1}
\newcommand{\de}{\delta}
\newcommand{\FF}{\mathbb{F}}
\newcommand{\G}{{\rm G}_2}
\newcommand{\ga}{\gamma}
\newcommand{\la}{\lambda}
\newcommand{\matr}[4]{\left(\begin{array}{cc}
#1 & #2 \\
#3 & #4 \\
\end{array}\right)}
\newcommand{\OO}{\mathbf{O}}
\newcommand{\ov}[1]{\overline{#1}}
\newcommand{\RR}{\mathbb{R}}
\newcommand{\uu}{\mathbf{u}}
\newcommand{\vv}{\mathbf{v}}
\newcommand{\zero}{\mathbf{0}}

\title{On polynomial equations over split-octonions: the arbitrary field case}

\authors{Artem Lopatin}

\authorinfo{State University of Campinas (UNICAMP), Brazil}{dr.artem.lopatin@gmail.com}

\abstract{%
    Over the split-octonion algebra defined over an arbitrary field, we solve all polynomial equations whose coefficients are scalar except for the constant term. As an application, we determine the square and cubic roots of an octonion.
    }



\VOLUME{34}
\YEAR{2026}
\ISSUE{2}
\NUMBER{12}
\DOI{https://doi.org/10.46298/cm.17309}
\licence{CC BY-SA 4.0}
\editinfo{January 15, 2026}{March 30, 2026}{Jaqueline Mesquita, Mariel Sáez, Rafael Potrie and Tiago Macedo}

\acknowledgments{These results were obtained during a visit to the Sino-Russian Mathematics Center at Peking University in the autumn of 2025. We gratefully acknowledge the support and hospitality provided during this visit.}

\begin{document}

\section{Introduction}

Unless stated otherwise, let $\FF$ be a field  (which may be finite) of arbitrary characteristic $\Char\FF \geq 0$. All vector spaces and algebras are assumed to be defined over~$\FF$.

\subsection{Equations over octonions}

The problem of solving polynomial equations has historically been one of the central problems in mathematics and has played a fundamental role in the development of algebraic geometry and other branches of mathematics. Polynomial equations have been studied not only over fields, but also over matrix algebras, quaternion algebras, octonion algebras, and other noncommutative or nonassociative structures.

In general, an \emph{octonion algebra} $\mathbf{C}$ (also called a \emph{Cayley algebra}) over a field $\FF$ is a nonassociative alternative unital algebra of dimension $8$, endowed with a non-singular quadratic multiplicative form
\[
n \colon \mathbf{C} \to \FF,
\]
called the \emph{norm}. The norm $n$ is said to be
\begin{enumerate}
\item[$\bullet$] \emph{isotropic} if $n(a) = 0$ for some non-zero $a \in \mathbf{C}$. In this case, there exists a unique octonion algebra $\OO_{\FF}$ over $\FF$ with isotropic norm (see Theorem~1.8.1 in~\cite{Springer_Veldkamp_book_2000}). This algebra is called the \emph{split-octonion algebra}.

\item[$\bullet$] \emph{anisotropic}, otherwise. In this case, the octonion algebra $\mathbf{C}$ is a division algebra.
\end{enumerate}

Note that if the field $\FF$ is algebraically closed, then any octonion algebra over $\FF$ is isomorphic to the split-octonion algebra $\OO_{\FF}$ (see, for example, Lemma~2.2 in~\cite{Lopatin_Zubkov_Eq_over_O}). By Artin's theorem, in any alternative algebra every subalgebra generated by two elements is associative. Consequently, every octonion algebra is power-associative; that is, the subalgebra generated by a single element is associative. Therefore, for any $a \in \mathbf{C}$, the power $a^n$ is well defined without the need to specify the placement of parentheses.

 Polynomial equations over octonion algebras have been studied in various settings. In particular, Rodríguez-Ordóñez~\cite{Rodriguez-Ordonez_2007} proved that every polynomial equation of positive degree over the algebra $\mathbf{A}_{\RR}$ of \emph{Cayley numbers} (i.e., the division algebra of real octonions), with the only highest-degree term present, has at least one solution. An explicit algorithm for solving quadratic equations of the form $x^2 + b x + c = 0$  
over $\mathbf{A}_{\RR}$ was provided by Wang, Zhang, and Zhang~\cite{Wang_Zhang_2014}, along with criteria determining whether such an equation has one, two, or infinitely many solutions.


Flaut and Shpakivskyi~\cite{Flaut_Shpakivskyi_2015} studied the equation $x^n = a$ over real octonion division algebras. For an octonion division algebra $\mathbf{C}$ over an arbitrary field $\FF$, Chapman~\cite{Chapman_2020_JAA} developed a complete method for finding the solutions of a general polynomial equation
$a_n x^n + a_{n-1} x^{n-1} + \cdots + a_1 x + a_0 = 0$ 
over $\mathbf{C}$. Furthermore, Chapman and Vishkautsan~\cite{Chapman_Vishkautsan_2022}, working over a division algebra $\mathbf{C}$, determined the solutions of the polynomial equation
$(a_n c) x^n + (a_{n-1} c) x^{n-1} + \cdots + (a_1 c) x + (a_0 c) = 0$, 
and also discussed the solutions of the related equation
$(c a_n) x^n + (c a_{n-1}) x^{n-1} + \cdots + (c a_1) x + (c a_0) = 0$.
Chapman and Levin~\cite{Chapman_Levin_2023} introduced a method for finding so-called \emph{alternating roots} of polynomials over an arbitrary division Cayley--Dickson algebra. In a subsequent work, Chapman and Vishkautsan~\cite{Chapman_Vishkautsan_2025} investigated conditions under which, for a root $a$ of a polynomial $f(x)$ over a general Cayley--Dickson algebra, there exists a factorization $f(x) = g(x)(x - a)$
for some polynomial $g(x)$.
Working over the split-octonion algebra over an algebraically closed field, Lopatin and Rybalov~\cite{Lopatin_Rybalov_2025} solved all polynomial equations in which all coefficients except the constant term are scalar. As a consequence, the $n$-th roots of a split-octonion were computed. In~\cite{Lopatin_Zubkov_Eq_over_O}, Lopatin and Zubkov studied the linear equations given by $a x = c$, $(a x) b = c$ and $a (b x) = c$ over the split-octonion algebra $\OO$ when the base field $\FF$ is algebraically closed. 
It is worth noting that, over a division octonion algebra, these equations are easily solvable and admit a unique solution whenever $a,b \neq 0$. In contrast, over an algebraically closed field the situation is substantially more delicate. As a consequence of the main result in~\cite{Lopatin_Zubkov_Eq_over_O}, it was shown that if a linear monomial equation over octonions with a non-zero constant term has at least two solutions, then it necessarily admits an invertible solution.


The split-octonion algebra has numerous applications in physics. For example, the Dirac equation, which describes the motion of a free spin-$\tfrac12$ particle such as an electron or a proton, can be formulated in terms of split-octonions (see~\cite{Koplinger_2006, Gogberashvili_2006, Gogberashvili_2006_Dirac, Koplinger_2007}). Further applications of split-octonions arise in electromagnetic theory (see~\cite{Chanyal_2017_CommTP, Chanyal_Bisht_Negi_2013, Chanyal_Bisht_Negi_2011}), geometrodynamics (see~\cite{Chanyal_2015_RepMP}), unified quantum theories (see~\cite{Krasnov_2022, Bisht_Dangwal_Negi_2008, Castro_2007}), and special relativity (see~\cite{Gogberashvili_Sakhelashvili_2015}).

Polynomial equations over arbitrary algebras have recently been investigated by Illmer and Netzer~\cite{Illmer_Netzer_2024}, who established conditions guaranteeing the existence of a common solution to $n$ polynomial equations in $n$ variables, with an application to polynomial equations over $\mathbf{A}_{\RR}$. Linear equations over matrix algebras have also been studied extensively (see, for example, \cite{Eq_BK_79, Eq_BK_80, Eq_B_99, Eq_T_00, Eq_L_06, Eq_FLLT_19}). The main questions addressed in these works include
\begin{enumerate}
\item[$\bullet$] determining conditions for the existence of solutions to linear equations;
\item[$\bullet$] describing the general form of the solutions.
\end{enumerate}

\subsection{Results}

In Section~\ref{section_O} we define the octonion algebra $\OO$ and its automorphism group $\Aut(\OO)$. In Section~\ref{section_main} we extend the results of~\cite{Lopatin_Rybalov_2025} from the case of an algebraically closed field to the case of an arbitrary field, which includes, in particular, the case of $\RR$ that is important for applications in physics.

More precisely, we solve the equation
\begin{equation}\label{eq_main}
\al_n x^n + \al_{n-1} x^{n-1} + \cdots + \al_1 x = c
\end{equation}%
\noindent
with scalar coefficients $\al_1,\ldots,\al_n \in \FF$ and a possibly non-scalar constant term $c \in \OO$, where the variable is $x \in \OO$ (see Theorem~\ref{theo_main}). The solution of equation~(\ref{eq_main}) is reduced to solving polynomial equations over the field $\FF$. In Corollary~\ref{cor_main} we provide further details for the case where $c \in \OO$ is non-scalar. As applications, we consider
\begin{enumerate}
\item[$\bullet$] the quadratic equation $x^2 = c$ in Proposition~\ref{prop_2};
\item[$\bullet$] the cubic equation $x^3 = c$ over the real numbers in Proposition~\ref{prop_3}.
\end{enumerate}
\noindent
In Corollary~\ref{cor_number_of_solutions} we show that, when $c$ is non-scalar, the number of solutions of equation~(\ref{eq_main}) is finite.

Besides allowing an arbitrary base field, the main difference from the results of~\cite{Lopatin_Rybalov_2025} is that we do not assume $c$ to be in a \emph{canonical} form. Consequently, our proofs differ substantially from those in~\cite{Lopatin_Rybalov_2025}. Moreover, when we employ a canonical form of an octonion (see Proposition~\ref{prop_O_canon}), it is different from the canonical forms considered in~\cite{Lopatin_Rybalov_2025} (see Definition~2.3 in~\cite{Lopatin_Rybalov_2025} for details).

\section{Octonions}\label{section_O}

\corr{The definitions presented in this section are taken from~\cite{LZ_2}. Additional material on octonions can be found in the books~\cite{Zhevlakov_book, McCrimmon_book}.}

\subsection{Split-octonions}\label{section_split}

The \emph{split octonion algebra} $\OO=\OO(\FF)$, also known as the \emph{split Cayley algebra}, is the vector space consisting of all matrices of the form
\[
a=\matr{\al}{\uu}{\vv}{\be},
\]
where $\al,\be \in \FF$ and $\uu,\vv \in \FF^{3}$, endowed with the multiplication
\[
a a' =
\matr{\al \al' + \uu \cdot \vv'}{\al \uu' + \be' \uu - \vv \times \vv'}{\al' \vv + \be \vv' + \uu \times \uu'}{\be \be' + \vv \cdot \uu'},
\qquad
a'=\matr{\al'}{\uu'}{\vv'}{\be'}.
\]
Here
\[
\uu \cdot \vv = u_1 v_1 + u_2 v_2 + u_3 v_3,
\qquad
\uu \times \vv = (u_2 v_3 - u_3 v_2,\; u_3 v_1 - u_1 v_3,\; u_1 v_2 - u_2 v_1).
\]

\noindent{}For brevity, we denote by $\cc_1=(1,0,0)$, $\cc_2=(0,1,0)$, $\cc_3=(0,0,1)$, and $\zero=(0,0,0)$ the standard basis vectors of $\FF^3$ and the zero vector, respectively. Consider the following basis of $\OO$:
\[
e_1=\matr{1}{\zero}{\zero}{0}, \;
e_2=\matr{0}{\zero}{\zero}{1}, \;
\uu_i=\matr{0}{\cc_i}{\zero}{0}, \;
\vv_i=\matr{0}{\zero}{\cc_i}{0},\;
i=1,2,3.
\]
The unity element of $\OO$ is given by $1_{\OO}=e_1+e_2$. Note that the multiplication in this basis satisfies
\[
\uu_i \uu_j = (-1)^{\epsilon_{ij}} \vv_k,
\qquad
\vv_i \vv_j = (-1)^{\epsilon_{ji}} \uu_k,
\]
where $\{i,j,k\}=\{1,2,3\}$ and $\epsilon_{ij}$ denotes the parity of the permutation
\[
\begin{pmatrix}
1 & 2 & 3 \\
k & i & j
\end{pmatrix}.
\]


The algebra $\OO$ is endowed with a linear involution
\[
\ov{a}=\matr{\be}{-\uu}{-\vv}{\al},
\]
which satisfies $\ov{a a'}=\ov{a'}\,\ov{a}$. The associated \emph{norm} is defined by
\[
n(a)=\al\be-\uu\cdot\vv,
\]
and it induces a nondegenerate symmetric bilinear \emph{form}
\[
q(a,a')=n(a+a')-n(a)-n(a')
= \al\be' + \al'\be - \uu\cdot \vv' - \uu'\cdot \vv .
\]
The linear \emph{trace} function is given by $\tr(a)=\al+\be$. Note that
\[
\tr(a)\,1_{\OO}=a+\ov{a}
\qquad \text{and} \qquad
n(a)\,1_{\OO}=a\,\ov{a}.
\]

\noindent{}The subspace of traceless octonions is denoted by
\[
\OO_0=\{a\in\OO \mid \tr(a)=0\},
\]
and the affine variety of octonions of zero norm by
\[
\OO_{\#}=\{a\in\OO \mid n(a)=0\}.
\]
The following identities hold for all $a,a'\in\OO$:
\begin{equation}\label{eq1}
\tr(a a')=\tr(a' a),
\qquad
n(a a')=n(a)\,n(a'),
\end{equation}
and
\begin{equation}\label{eq_nab}
n(a+a')=n(a)+n(a')-\tr(a a')+\tr(a)\tr(a').
\end{equation}

\noindent{}Moreover, every element $a\in\OO$ satisfies the quadratic identity
\begin{equation}\label{eq2}
a^2-\tr(a)\,a+n(a)\,1_{\OO}=0.
\end{equation}

The algebra $\OO$ is a simple \emph{alternative} algebra; that is, for all $a,b\in\OO$,
\begin{equation}\label{eq4}
a(ab)=(aa)b,
\qquad
(ba)a=b(aa).
\end{equation}
Furthermore, the involution interacts with multiplication via
\begin{equation}\label{eq4b}
\ov{a}(ab)=n(a)\,b,
\qquad
(ba)\ov{a}=n(a)\,b.
\end{equation}

The following remark is well known and can be proved easily.

\begin{remark}\label{remark_inv}
Let $a \in \OO$. Then exactly one of the following cases occurs:
\begin{enumerate}
\item[$\bullet$] If $n(a)\neq 0$, then there exist unique elements $b,c \in \OO$ such that
$b a = 1_{\OO}$ and
$a c = 1_{\OO}$.
In this case, $b = c = \ov{a}/n(a)$, 
and we denote this element by $a^{-1}$.

\item[$\bullet$] If $n(a)=0$, then $a$ admits neither a left inverse nor a right inverse in $\OO$.
\end{enumerate}
\end{remark}

As a direct consequence of identities~(\ref{eq4b}), for every $a \in \OO \setminus \OO_{\#}$ we have
\begin{equation}\label{eq_inv}
a^{-1}(a b) = b,
\qquad
(b a) a^{-1} = b,
\end{equation}
for all $b \in \OO$.

\subsection{The group $\Aut(\OO)$}\label{section_G2}

The group $\Aut(\OO)$ of all automorphisms of the algebra $\OO$ is the exceptional simple algebraic group $\G(\FF)$ when the field $\FF$ is algebraically closed. In the general case, the group $\Aut(\OO)$ contains a subgroup isomorphic to $\SL_3(\FF)$. More precisely, every element $g \in \SL_3(\FF)$ defines an automorphism of $\OO$ by
\[
a \longmapsto \matr{\al}{\uu g}{\vv g^{-T}}{\be},
\]
where $g^{-T}$ denotes $(g^{-1})^T$ and the vectors $\uu,\vv \in \FF^3$ are regarded as row vectors. For every $\uu,\vv \in \FF^3$, define automorphisms $\de_1(\uu), \de_2(\vv) \in \Aut(\OO)$ by
\[
\de_1(\uu)(a')
=
\matr{\al' - \uu \cdot \vv'}
{(\al' - \be' - \uu \cdot \vv')\uu + \uu'}
{\vv' - \uu' \times \uu}
{\be' + \uu \cdot \vv'},
\]
\[
\de_2(\vv)(a')
=
\matr{\al' + \uu' \cdot \vv}
{\uu' + \vv' \times \vv}
{(-\al' + \be' - \uu' \cdot \vv)\vv + \vv'}
{\be' - \uu' \cdot \vv},
\]
where $a'=\matr{\al'}{\uu'}{\vv'}{\be'}$.

A direct computation shows that the map
\begin{equation}\label{eq_h}
\hbar \colon \OO \to \OO,
\qquad
a \longmapsto \matr{\be}{-\vv}{-\uu}{\al},
\end{equation}
also belongs to $\Aut(\OO)$.

The action of $\Aut(\OO)$ on $\OO$ satisfies the following properties:
\begin{equation}\label{eq_properties}
\ov{g a} = g \ov{a},
\qquad
\tr(g a) = \tr(a),
\qquad
n(g a) = n(a),
\qquad
q(g a, g a') = q(a, a'),
\end{equation}
for all $g \in \Aut(\OO)$ and $a,a' \in \OO$ (see, for example, equation~(2.2) in~\cite{Elduque_2023}). Consequently, the group $\Aut(\OO)$ also acts naturally on the subspaces $\OO_0$ and $\OO_{\#}$.

\subsection{Canonical octonions}\label{section_canon}

A minimal set of representatives for the $\Aut(\OO)$-orbits in $\OO$ was described in Proposition~3.3 of~\cite{Lopatin_Zubkov_2025} in the case where the base field is algebraically closed. We extend this result to an arbitrary field in the following proposition. The elements appearing in Proposition~\ref{prop_O_canon} will be referred to as \emph{canonical octonions}.

\begin{proposition}\label{prop_O_canon}
A minimal set of representatives for the $\Aut(\OO)$-orbits in $\OO$ consists of the following elements:
\begin{enumerate}
\item[1.] $\al\,1_{\OO}$,

\item[2.] $\displaystyle \matr{\al}{(\be,0,0)}{(1,0,0)}{0} \in \OO$,
\end{enumerate}
where $\al,\be \in \FF$. In other words, $\OO$ is the disjoint union of the following $\Aut(\OO)$-orbits:
\[
\al\,1_{\OO}
\qquad \text{and} \qquad
O(\al,\be) := \Aut(\OO)\cdot \matr{\al}{(\be,0,0)}{(1,0,0)}{0},
\]
with $\al,\be \in \FF$.
\end{proposition}

\begin{proof}
By Lemma~3.3 of~\cite{Elduque_2023}, two elements $a,b \in \OO \setminus \FF 1_{\OO}$ belong to the same $\Aut(\OO)$-orbit if and only if $\tr(a)=\tr(b)$ and $n(a)=n(b)$. The proof is completed by observing that
\[
\tr\!\left(\matr{\al}{(\be,0,0)}{(1,0,0)}{0}\right)=\al,
\qquad
n\!\left(\matr{\al}{(\be,0,0)}{(1,0,0)}{0}\right)=-\be.
\]
\end{proof}

The above proof of Proposition~\ref{prop_O_canon} is non-constructive and relies on
Lemma~3.3 of~\cite{Elduque_2023}. For the sake of completeness, we present an
alternative proof that does not depend on Lemma~3.3 of~\cite{Elduque_2023} and
that, for every $a\in\OO$, explicitly constructs an element
$g \in \Aut(\OO)$ such that $ga$ is in canonical form. This constructive
approach may be useful for practical computations involving octonions.

\begin{proof}[Constructive proof of Proposition~\ref{prop_O_canon}]
In this proof, we use the symbol \( \ast \) to denote an arbitrary element of
\( \FF \). Let
$a=\matr{\al_1}{\uu}{\vv}{\al_8}$.
Then one of the following two cases occurs.

\medskip
\noindent\textbf{1.} Assume that $\uu=\vv=\zero$.
If $\al_1=\al_8$, then $a\in \FF 1_{\OO}$ is canonical.

Now assume that $\al_1\neq \al_8$.
By acting on $a$ with $\de_2(1,0,0)$, we may assume that
$a=\matr{\al_1}{\zero}{(\al_8-\al_1,0,0)}{\al_8}$.
Applying part~(b) of Lemma~3.2 from~\cite{Lopatin_Zubkov_2025}, we may further
assume that
$a=\matr{\al_1}{\zero}{(1,0,0)}{\al_8}$.

If $\al_8=0$, then $a$ is canonical.
Otherwise, acting on $a$ with $\de_1(-\al_8,0,0)$, we obtain the canonical
octonion
$\matr{\al_1+\al_8}{(-\al_1\al_8,0,0)}{(1,0,0)}{\corr{0}}$.

\medskip

\noindent\textbf{2.} Assume that $\uu$ or $\vv$ is non-zero. Applying the
automorphism $\hbar$, we may assume that $\vv$ is non-zero. Using part~(a) of
Lemma~3.2 from~\cite{Lopatin_Zubkov_2025} together with the automorphism $\hbar$,
we reduce to the case $\vv=(1,0,0)$. Acting by $\de_1(-\al_8,0,0)$, we then obtain
the octonion
$$
\matr{\al_1+\al_8}{\uu'}{(1,\ast,\ast)}{0}.
$$
As above, similarly to part~(a) of Lemma~3.2 from~\cite{Lopatin_Zubkov_2025}, we further reduce to the octonion
$$
\matr{\al}{\uu''}{(1,0,0)}{0},
$$
where $\al=\al_1+\al_8$ and $\uu''$ is one of the following two vectors:
\begin{enumerate}
\item[$\bullet$] $\uu''=(\ast,0,0)$. In this case, we have obtained a canonical
octonion.

\item[$\bullet$] $\uu''=(0,1,0)$. If $\al\neq 0$, we apply
$\de_1\!\left(0,-\frac{1}{\al},0\right)$ to obtain the canonical octonion
$\matr{\al}{\zero}{(1,0,0)}{0}$. If $\al=0$, we apply $\de_2(0,0,1)$ to obtain
the canonical octonion $\matr{0}{\zero}{(1,0,0)}{0}$,
\end{enumerate}
and the proof is finished.
\end{proof}

\section{Polynomial equations}\label{section_main}

Recursively, we define the (commutative and associative) Generalized Fibonacci polynomials 
$p_n=p_n(y,z)\in\FF[y,z]$ for all $n\geqslant -1$ by
$$
p_{-1}=0,\qquad p_0=1,\qquad p_{k+1}=y p_k+ z p_{k-1}
\quad \text{for all } k\geqslant 0.
$$
In particular,
$$
\begin{array}{rcl}
p_1 &=& y, \\
p_2 &=& y^2+z,\\ 
p_3 &=& y^3+2yz, \\ 
p_4 &=& y^4 + 3 y^2 z + z^2, \\ 
p_5 &=& y^5+4 y^3 z + 3 y z^2,\\ 
p_6 &=& y^6+5 y^4 z + 6 y^2 z^2 + z^3. \\  
\end{array}
$$

Let $f(y)=\al_n y^n + \cdots + \al_1 y\in\FF[y]$ be a non-zero polynomial without
constant term, where $\al_1,\ldots,\al_n\in\FF$, $\al_n\neq 0$, and $n\geqslant 1$.
We define the polynomials $\widehat{f}(y,z)$ and $\widecheck{f}(y,z)$ in
$\FF[y,z]$ by
$$
\begin{array}{rcl}
\widehat{f}(y,z) &=& \al_n p_{n-1}(y,z) + \cdots + \al_2 p_1(y,z) + \al_1 p_0(y,z), \\
\widecheck{f}(y,z) &=& \al_n p_{n-2}(y,z) + \cdots + \al_2 p_0(y,z) + \al_1 p_{-1}(y,z).
\end{array}
$$

For each $x\in\OO$, we naturally define the substitution
$$
f(x)=\al_n x^n + \cdots + \al_1 x \in \OO.
$$

\begin{proposition}\label{prop_power}
Given $a\in\OO$, we denote $\tr(a)=\al$ and $n(a)=-\be$, where $\al,\be\in\FF$.
\begin{enumerate}
\item[1.] For $n\geqslant 1$, we have
\begin{equation}\label{eq_power}
a^n = p_{n-1}(\al,\be)\, a + p_{n-2}(\al,\be)\,\be 1_{\OO}.
\end{equation}

\item[2.] For a non-zero polynomial $f(y)\in\FF[y]$ without constant term, we
have
\begin{equation}\label{eq_fa}
f(a) = \widehat{f}(\al,\be)\, a + \widecheck{f}(\al,\be)\,\be 1_{\OO}.
\end{equation}
\end{enumerate}
\end{proposition}

\begin{proof}
\textbf{1.} For brevity, set $q_n=p_n(\al,\be)$. The quadratic equation~(\ref{eq2})
can be rewritten in the form
\begin{equation}\label{eq2_new}
a^2 = \al a + \be 1_{\OO}.
\end{equation}

We prove formula~(\ref{eq_power}) by induction on $n\geqslant 1$. The case $n=1$
is trivial.

Assume that~(\ref{eq_power}) holds for $n=k$. Then, by the induction hypothesis
and~(\ref{eq2_new}), we obtain
$$
a^{k+1}
= q_{k-1} a^2 + q_{k-2}\be a
= \big(q_{k-1}\al + q_{k-2}\be\big)a + \be q_{k-1} 1_{\OO}.
$$
By the recursive definition of $p_n(y,z)$, this proves part~1.

\medskip
\noindent\textbf{2.}
Let $f(y)=\sum_{k=1}^n \al_k y^k$, where $\al_k\in\FF$. Applying part~1, we obtain
$$
\begin{aligned}
f(a)
&= \sum_{k=1}^n \al_k p_{k-1}(\al,\be)\, a
 + \sum_{k=1}^n \al_k p_{k-2}(\al,\be)\,\be 1_{\OO} \\
&= \widehat{f}(\al,\be)\, a + \widecheck{f}(\al,\be)\,\be 1_{\OO},
\end{aligned}
$$
which proves~(\ref{eq_fa}).
\end{proof}

Proposition~\ref{prop_power} implies the following remark.

\begin{remark}\label{rem_power}
Let 
$a=\matr{\al}{(\be,0,0)}{(1,0,0)}{0}\in\OO$
for some $\al,\be\in\FF$, and let $n\geqslant 2$. Then
$$
a^n = \matr{
p_n(\al,\be)
}{
(\be\, p_{n-1}(\al,\be),0,0)
}{
(p_{n-1}(\al,\be),0,0)
}{
\be\, p_{n-2}(\al,\be)
}.
$$
\end{remark}

\corr{
Recall that $O(\al,\be)$ is defined in Proposition~\ref{prop_O_canon} as the $\Aut(\OO)$-orbit of the canonical octonion $\matr{\al}{(\be,0,0)}{(1,0,0)}{0}$.
}

\begin{theorem}\label{theo_main}
Assume that $f(y)\in\FF[y]$ is a non-zero polynomial without constant term, 
and let $c\in\OO$. Let $X\subseteq \OO$ be the set of all solutions of the 
equation $f(x)=c$, where $x\in\OO$ is a variable.

\begin{enumerate}
\item[1.] Assume $c=\ga 1_{\OO}$ for some $\ga\in\FF$. Then
$$
X =  \big\{\nu 1_{\OO} \,\big|\, \nu\in\FF \text{ satisfies } f(\nu)=\ga\big\} 
\;\cup 
\bigcup_{\substack{\la,\mu \in \FF \\ \widehat{f}(\la,\mu)=0 \\ \mu\,\widecheck{f}(\la,\mu)=\ga}} 
O(\la,\mu).
$$


\item[2.] Assume $c\not\in \FF 1_{\OO}$. Then $x\in X$ if and only if there exist 
$\la,\mu\in\FF$ satisfying
\begin{equation}\label{eq_system}
\left\{
\begin{array}{rcl} 
\la \,\widehat{f}(\la,\mu)  + 2 \mu \,\widecheck{f}(\la,\mu)  & = & \tr(c),\\[1mm]
-\mu \,\widehat{f}(\la,\mu)^2 + \la \mu \,\widehat{f}(\la,\mu) \widecheck{f}(\la,\mu) 
+ \mu^2 \,\widecheck{f}(\la,\mu)^2  & = & n(c),\\[1mm]
\widehat{f}(\la,\mu)  & \neq & 0,
\end{array}
\right.
\end{equation}
such that
\begin{equation}\label{eq_solution} 
x = \frac{1}{\widehat{f}(\la,\mu)}\Big( c - \mu \,\widecheck{f}(\la,\mu) 1_{\OO} \Big).
\end{equation}
In this case, we have $\tr(x)=\la$ and $n(x)=-\mu$.
\end{enumerate}
\end{theorem}
\begin{proof}
\noindent\textbf{1.} Let $c=\ga 1_{\OO}$ for some $\ga\in\FF$, and consider an
arbitrary $x\in \OO$. By Proposition~\ref{prop_O_canon}, there exists
$g\in\Aut(\OO)$ such that $gx$ is canonical. Note that the equality $f(x)=c$
is equivalent to
\begin{equation}\label{eq_gxc1}
f(gx)= \ga 1_{\OO}.
\end{equation}

\noindent One of the following possibilities occurs:

\begin{enumerate}
\item[{\bf(a)}] $gx=\nu 1_{\OO}$ for some $\nu\in \FF$. In this case,
equality~(\ref{eq_gxc1}) is equivalent to $f(\nu)=\ga$. Hence,
$$
x = g^{-1}(\nu 1_{\OO}) = \nu 1_{\OO}.
$$

\item[{\bf(b)}] $gx=\matr{\la}{(\mu,0,0)}{(1,0,0)}{0}$ for some
$\la,\mu\in \FF$. By part~2 of Proposition~\ref{prop_power}, equality~(\ref{eq_gxc1}) is equivalent to
$$
\left\{
\begin{array}{rcl} 
\widehat{f}(\la,\mu) &= &0,\\ [1mm]
\mu\,\widecheck{f}(\la,\mu) &= &\ga. \\
\end{array}\right.
$$
Since $x\in O(\la,\mu)$, we have that $x\in X$ if and only if
$O(\la,\mu)\subseteq X$.
\end{enumerate}

\noindent The claim is thus proven.

\medskip
\noindent\textbf{2.} Assume $c\not\in \FF 1_{\OO}$.

\medskip
\noindent\textbf{(a)} Assume that $x\in X$. By part~2 of Proposition~\ref{prop_power}, the equality $f(x)=c$ is equivalent to
\begin{equation}\label{eq_main2}
\widehat{f}(\la,\mu)\, x + \mu\, \widecheck{f}(\la,\mu)\, 1_{\OO} = c,
\end{equation}%
\noindent{}where $\la = \tr(x)$ and $\mu = - n(x)$. If $\widehat{f}(\la,\mu)=0$, then
equality~(\ref{eq_main2}) would imply $c\in \FF 1_{\OO}$, a contradiction.
Hence, $x$ is given by equality~(\ref{eq_solution}).  

Applying the trace and norm to equality~(\ref{eq_main2}) and using the linearity of
the trace together with formula~(\ref{eq_nab}), we obtain that system~(\ref{eq_system}) holds.

\medskip
\noindent\textbf{(b)} Conversely, assume that $x\in\OO$ is defined by
equality~(\ref{eq_solution}) for some $\la,\mu\in\FF$ satisfying
system~(\ref{eq_system}).  System~(\ref{eq_system}) implies
$$
\tr(x) = \frac{1}{\widehat{f}(\la,\mu)} \Big( \tr(c) - 2 \mu\, \widecheck{f}(\la,\mu) \Big) = \la.
$$
Similarly, formula~(\ref{eq_nab}) together with system~(\ref{eq_system}) gives
\begin{multline*}
n(x) = \frac{1}{\widehat{f}(\la,\mu)^2} \Big( n(c) + \mu^2 \widecheck{f}(\la,\mu)^2 
- \mu\, \tr(c) \widecheck{f}(\la,\mu) \Big) \\
= \frac{1}{\widehat{f}(\la,\mu)^2} \Big( n(c) 
- \la \mu\, \widehat{f}(\la,\mu) \widecheck{f}(\la,\mu)
- \mu^2 \widecheck{f}(\la,\mu)^2 \Big) = -\mu.
\end{multline*}

\noindent{}Therefore, by part~2 of Proposition~\ref{prop_power} and equality~(\ref{eq_solution}),
$$
f(x) = \widehat{f}(\la,\mu)\, x + \mu\, \widecheck{f}(\la,\mu)\, 1_{\OO} = c.
$$
Hence, $x$ belongs to $X$.
\end{proof}

The following corollary is a straightforward consequence of Theorem~\ref{theo_main}. 
In its formulation and hereafter, we use the standard notation 
$\FF^{\times}$ to denote the set of all non-zero elements of $\FF$.

\begin{corollary}\label{cor_main} 
Assume that $f(y)\in\FF[y]$ is a non-zero polynomial without constant term, and 
let $c\in\OO \setminus \FF 1_{\OO}$. Let $X\subseteq \OO$ be the set of all 
solutions of the equation $f(x)=c$, where $x\in\OO$ is a variable.

\begin{enumerate}
\item[1.] Assume $\Char\FF\neq 2$. Then $x\in X$ if and only if there exist 
$\la,\mu\in\FF$ satisfying
\begin{equation}\label{eq_cor_sys1}
\left\{
\begin{array}{rcl} 
(\la^2 + 4\mu)  \widehat{f}(\la,\mu)^2 & = & \tr(c)^2 - 4\,n(c),\\[1mm]
\la \widehat{f}(\la,\mu) + 2 \mu \widecheck{f}(\la,\mu)  & = & \tr(c),\\[1mm]
\widehat{f}(\la,\mu)  & \neq & 0,
\end{array}
\right.
\end{equation}
such that
\begin{equation}\label{eq_cor_solution1} 
x = \frac{1}{\widehat{f}(\la,\mu)}\Big( c - \frac{\tr(c)}{2} 1_{\OO} \Big) 
+ \frac{\la}{2}\, 1_{\OO}.
\end{equation}

\item[2.] Assume $\Char\FF = 2$. 

\begin{enumerate}
\item[(a)] If $\tr(c)\neq0$, then $x\in X$ if and only if there exist 
$\la\in \FF^{\times}$ and $\mu\in\FF$ satisfying
\begin{equation}\label{eq_cor_sys2a}
\left\{
\begin{array}{rcl} 
\mu^2  \widecheck{f}(\la,\mu)^2 + \mu\,\tr(c)\,\widecheck{f}(\la,\mu) & = & n(c) + \frac{\mu}{\la^2} \tr(c)^2,\\[1mm]
\la \widehat{f}(\la,\mu)  & = & \tr(c),
\end{array}
\right.
\end{equation}
such that
\begin{equation}\label{eq_cor_solution2a} 
x = \frac{\la}{\tr(c)} \Big( c + \mu\, \widecheck{f}(\la,\mu) 1_{\OO} \Big).
\end{equation}

\item[(b)] If $\tr(c)=0$, then $x\in X$ if and only if there exists $\mu\in\FF$ satisfying
\begin{equation}\label{eq_cor_sys2b}
\left\{
\begin{array}{rcl} 
\mu \widehat{f}(0,\mu)^2 + \mu^2 \widecheck{f}(0,\mu)^2 & = & n(c),\\[1mm]
\widehat{f}(0,\mu)  & \neq & 0,
\end{array}
\right.
\end{equation}
such that
\begin{equation}\label{eq_cor_solution2b} 
x = \frac{1}{\widehat{f}(0,\mu)} \Big( c + \mu\, \widecheck{f}(0,\mu) 1_{\OO} \Big).
\end{equation}
\end{enumerate}
\end{enumerate}

\noindent Note that in each case we have $\tr(x)=\la$ and $n(x)=-\mu$, 
where in case~2(b) we have $\la=0$.
\end{corollary}


\begin{corollary}\label{cor_number_of_solutions} 
Assume that $f(y)\in\FF[y]$ is a non-zero polynomial of degree $n\geqslant 1$ 
without constant term, and let $c\in\OO\setminus \FF 1_{\OO}$. 
Denote by $X\subseteq \OO$ the set of all solutions of the equation $f(x)=c$. 
Then
$$
|X| \leqslant n^2.
$$
\end{corollary}
\begin{proof}
Consider the inclusion $\FF\subseteq \overline{\FF}$, where $\overline{\FF}$ is 
the algebraic closure of $\FF$. Let $X'\subseteq \OO(\overline{\FF})$ denote 
the set of all solutions of $f(x)=c$ in $\OO(\overline{\FF})$. 

Since $c\not\in \FF 1_{\OO}$, Corollary~3.3 of~\cite{Lopatin_Rybalov_2025} 
implies that $|X'|\leqslant n^2$. 

Finally, since $\OO\subseteq \OO(\overline{\FF})$, we have $X\subseteq X'$, 
which proves the claim.
\end{proof}

\section{Roots of octonions}\label{section_root}

\begin{definition}
\corr{
\begin{enumerate}
\item[$\bullet$] We denote by $\FF^2$ the set $\{\al^2 \mid \al \in \FF\}$.

\item[$\bullet$] Let $\be \in \FF^2$. We fix an element $\sqrt{\be}$ from the non-empty set 
$\{\al \in \FF \mid \al^2 = \be\}$. 
Observe that the set of all solutions to the equation $\nu^2 = \be$, where $\nu \in \FF$ is a variable, is $\{\pm \sqrt{\be}\}$.
\end{enumerate}
}
\end{definition}

\begin{proposition}\label{prop_2}
Let $X\subseteq\OO$ be the set of all solutions of the equation $x^2=c$, 
where $x\in\OO$ is a variable and $c\in\OO$.

\begin{enumerate}
\item[1.] Assume $c = \ga 1_{\OO}$ for some $\ga\in\FF$. Then
$$
X = \corr{\pm \sqrt{\ga}}1_{\OO} \;\cup\; O(0,\ga).
$$

\item[2.] Assume $c\not\in \FF 1_{\OO}$ and $\Char\FF \neq 2$.

\begin{enumerate}
\item[(a)] If $n(c)\in \FF^2$, denote
$$
\al = \tr(c) + 2\sqrt{n(c)}, \qquad \be = \tr(c) - 2\sqrt{n(c)}.
$$

\begin{itemize}
\item If $\al,\be \in (\FF^{\times})^2$, then 
$$
X = \Big\{ \pm \frac{1}{\sqrt{\al}} \big(c+\sqrt{n(c)} 1_{\OO}\big), \;\pm 
            \frac{1}{\sqrt{\be}} \big(c-\sqrt{n(c)} 1_{\OO}\big) \Big\}.
$$

\item If $\al\in (\FF^{\times})^2$ and $\be \not\in (\FF^{\times})^2$, then
$$
X = \Big\{ \pm \frac{1}{\sqrt{\al}} \big(c+\sqrt{n(c)} 1_{\OO}\big) \Big\}.
$$

\item If $\al\not\in (\FF^{\times})^2$ and $\be \in (\FF^{\times})^2$, then
$$
X = \Big\{ \pm \frac{1}{\sqrt{\be}} \big(c-\sqrt{n(c)} 1_{\OO}\big) \Big\}.
$$

\item If $\al,\be \not\in (\FF^{\times})^2$, then $X = \emptyset$.
\end{itemize}

\item[(b)] If $n(c)\not\in \FF^2$, then $X = \emptyset$.
\end{enumerate}

\item[3.] Assume $c\not\in \FF 1_{\OO}$ and $\Char\FF = 2$.

\begin{enumerate}
\item[(a)] If $\tr(c)\neq 0$, then
$$
X = 
\begin{cases} 
\frac{1}{\sqrt{\tr(c)}} \Big(c + \sqrt{n(c)} 1_{\OO}\Big), & \text{if } 
\tr(c), n(c) \in \FF^2,\\[1mm]
\emptyset, & \text{otherwise.}
\end{cases}
$$

\item[(b)] If $\tr(c) = 0$, then $X = \emptyset$.
\end{enumerate}

\end{enumerate}
\end{proposition}
\begin{proof} 
Note that for $f(y)=y^2\in \FF[y]$ and every $\la,\mu\in\FF$, we have
$$
\widehat{f}(\la,\mu) = p_1(\la,\mu) = \la, \qquad 
\widecheck{f}(\la,\mu) = p_0(\la,\mu) = 1.
$$

\medskip
\noindent\textbf{1.} Assume $c=\ga 1_{\OO}$ for some $\ga\in\FF$. 
Then the required statement follows immediately from Theorem~\ref{theo_main}.

\medskip
\noindent\textbf{2.} Assume $c\not\in \FF 1_{\OO}$ and $\Char\FF\neq 2$. 
By Corollary~\ref{cor_main}, $x\in X$ if and only if 
\begin{equation*}
x = \frac{1}{\la}\Big( c - \frac{\tr(c)}{2} 1_{\OO}\Big) + \frac{\la}{2} 1_{\OO}
\end{equation*}
for some $\la\in\FF^{\times}$ and $\mu\in\FF$ satisfying
\begin{equation}\label{eq_cor_sys1_sqrt}
\left\{\begin{array}{rcl}
(\la^2 + 4\mu)\, \la^2 &=& \tr(c)^2 - 4\,n(c),\\
 \la^2+2\mu &=& \tr(c).\\
\end{array}\right.
\end{equation}
System~(\ref{eq_cor_sys1_sqrt}) is equivalent to
\begin{equation*}
\left\{\begin{array}{rcl}
(\la^2 - \tr(c))^2 &=& 4\, n(c),\\
\mu &=& \frac{1}{2} \big(\tr(c) - \la^2 \big),\\
\end{array}\right.
\end{equation*}
and the required statement follows.

\medskip
\noindent\textbf{3.} Assume $c\not\in \FF 1_{\OO}$ and $\Char\FF = 2$.  

\noindent\textbf{(a)} Let $\tr(c)\neq 0$. By Corollary~\ref{cor_main}, 
$x\in X$ if and only if
\begin{equation*}
x = \frac{\corr{\la}}{\tr(c)} \Big( c + \mu 1_{\OO} \Big)
\end{equation*}
for some $\la\in\FF^{\times}$ and $\mu\in\FF$ satisfying
\begin{equation}\label{eq_cor_sys2a_sqrt}
\left\{\begin{array}{rcl}
\mu^2 + \mu\, \tr(c) &=& n(c) + \frac{\mu}{\la^2}\, \tr(c)^2,\\
\la^2 &=& \tr(c).\\
\end{array}\right.
\end{equation}
System~(\ref{eq_cor_sys2a_sqrt}) is equivalent to
\begin{equation*}
\la^2 = \tr(c), \qquad \mu^2 = n(c),
\end{equation*}
and the required statement follows.

\medskip
\noindent\textbf{(b)} Let $\tr(c) = 0$. Assume that $X$ is non-empty. 
Then, for any $x\in X$, Corollary~\ref{cor_main} implies the existence of 
$\mu\in\FF$ such that $\widehat{f}(0,\mu)\neq 0$. But 
$\widehat{f}(0,\mu)=0$, which is a contradiction. Hence $X=\emptyset$.
\end{proof}

Proposition~\ref{prop_2} implies the following corollary.

\begin{corollary}\label{cor_prop_2}
Assume $\FF=\RR$ and let $c\in\OO$.  
Then the equation $x^2=c$ has no solutions in $\OO$ if and only if
$c\notin \FF 1_{\OO}$ and one of the following conditions holds:
\begin{enumerate}
\item[$\bullet$] $0 \leqslant 4n(c) \leqslant \tr(c)^2$ and $\tr(c)\leqslant 0$;

\item[$\bullet$] $n(c)<0$.
\end{enumerate}
\end{corollary}

\begin{proposition}\label{prop_3}
Assume $\FF=\RR$ and let $c\in\OO$.  
Let $X\subseteq\OO$ be the set of all solutions of the equation $x^3=c$,
where $x\in\OO$ is a variable.
\begin{enumerate}
\item[1.] Assume $c=\ga 1_{\OO}$ for some $\ga\in\RR$. Then
\[
X= \big\{\sqrt[3]{\ga}\,1_{\OO}\big\} \,\cup\, O\big(-\sqrt[3]{\ga},-\sqrt[3]{\ga^2}\big).
\]

\item[2.] Assume $c\notin \RR 1_{\OO}$ and $\tr(c)\neq0$.  
Then $x\in X$ if and only if there exists $\la\in\RR$ satisfying
\begin{equation}\label{eq_prop_3_2}
\left\{
\begin{array}{rcl}
\big(2\la^3 + \tr(c)\big)^2\big(\la^3 - 4\tr(c)\big)
& = &
27\la^3\big(4n(c) - \tr(c)^2\big),\\[2mm]
2\la^3 + \tr(c) & \neq & 0,
\end{array}
\right.
\end{equation}
such that
\begin{equation}\label{eq_prop3_solution2}
x
=
\frac{3\la}{2\la^3 + \tr(c)}
\bigg( c - \frac{\tr(c)}{2}1_{\OO}\bigg)
+ \frac{\la}{2}1_{\OO}.
\end{equation}

\item[3.] Assume $c\notin \RR 1_{\OO}$ and $\tr(c)=0$.
\begin{enumerate}
\item[(a)] If $n(c)>0$, then $X=\{x_1,x_2\}$.

\item[(b)] If $n(c)=0$, then $X=\emptyset$.

\item[(c)] If $n(c)<0$, then $X=\{x_1\}$,
\end{enumerate}
where
\[
x_1 = -\frac{c}{\sqrt[3]{n(c)}},
\qquad
x_2 = \frac{1}{2}\bigg(
\frac{c}{\sqrt[3]{n(c)}} + \sqrt{3}\,\sqrt[6]{n(c)}\,1_{\OO}
\bigg).
\]
\end{enumerate}
\end{proposition}

\begin{proof}
Note that for $f(y)=y^3\in\RR[y]$ and every $\la,\mu\in\RR$, we have
\[
\widehat{f}(\la,\mu)=p_2(\la,\mu)=\la^2+\mu,
\qquad
\widecheck{f}(\la,\mu)=p_1(\la,\mu)=\la.
\]

\medskip
\noindent\textbf{I.}
Assume $c=\ga 1_{\OO}$ for some $\ga\in\RR$.
By Theorem~\ref{theo_main},
\[
X=
\big\{\sqrt[3]{\ga}\,1_{\OO}\big\}
\;\cup\;
\big\{O(\la,\mu)\mid \la,\mu\in\RR
\text{ satisfy }
\la^2+\mu=0
\text{ and }
\la\mu=\ga
\big\}.
\]
Considering separately the cases $\ga=0$ and $\ga\neq0$,
we obtain the statement of part~1.

\medskip
\noindent\textbf{II.}
Assume $c\notin \RR 1_{\OO}$.
By Corollary~\ref{cor_main}, an element $x\in\OO$ belongs to $X$ if and only if
\begin{equation}\label{eq_x}
x
=
\frac{1}{\la^2+\mu}
\bigg(c-\frac{\tr(c)}{2}1_{\OO}\bigg)
+
\frac{\la}{2}1_{\OO},
\end{equation}
for some $\la,\mu\in\RR$ satisfying
\begin{equation}\label{eq_sqrt3_sysII}
\left\{
\begin{array}{rcl}
(\la^2+4\mu)(\la^2+\mu)^2 &=& \tr(c)^2-4n(c),\\[2mm]
2\la\mu+\la(\la^2+\mu) &=& \tr(c),\\[2mm]
\la^2+\mu &\neq& 0.
\end{array}
\right.
\end{equation}
The second equation of system~(\ref{eq_sqrt3_sysII}) is equivalent to
\begin{equation}\label{eq_lamu}
3\la\mu=\tr(c)-\la^3.
\end{equation}

\medskip
\noindent\textbf{(a)}
Assume $\tr(c)\neq0$.
Then $\la\neq0$, since otherwise~(\ref{eq_lamu}) would imply $\tr(c)=0$,
a contradiction.
Hence
\[
\mu=\frac{\tr(c)-\la^3}{3\la}.
\]
Substituting this expression into system~(\ref{eq_sqrt3_sysII}) and equality~(\ref{eq_x}) 
yields system~(\ref{eq_prop_3_2}) and equality~(\ref{eq_prop3_solution2}), respectively. 
Note that the first equation of system~(\ref{eq_prop_3_2}) implies that $\la\neq0$. 
This completes the proof of part~2.

\medskip
\noindent\textbf{(b)}
Assume $\tr(c)=0$.
Considering separately the cases $\la=0$ and $\la\neq0$,
we obtain the statements of part~3.
\end{proof}

\end{document}